\documentclass{elsart}

\journal{Advances in Applied Mathematics}

\usepackage[utf8]{inputenc}
\usepackage{amsmath}
\usepackage{amssymb}
\usepackage[T1]{fontenc}
\usepackage[dvips]{graphicx}

\newcommand{\dfn}[1]{\textit{#1}} 
\newcommand{\ex}[1]{e(#1)} 
\newcommand{\vd}[1]{d(#1)} 

\begin{document}
\begin{frontmatter}

\title{Geometric representation of binary codes and computation of weight enumerators}

\author{Pavel Rytíř}
\ead{rytir@kam.mff.cuni.cz}
\address{Department of Applied Mathematics, Charles University in Prague, Malostranské nám. 25, Prague 118 00, Czech Republic}
\thanks{Supported by the Czech Science Foundation under the contract no. 201/09/H057}
\date{5 July 2009}


\begin{abstract}
For every linear binary code $\mathcal{C}$, we construct a geometric triangular configuration $\Delta$ so that the weight enumerator of $\mathcal{C}$ is obtained by a simple formula from the weight enumerator of the cycle space of $\Delta$. The triangular configuration $\Delta$ thus provides a geometric representation of $\mathcal{C}$ which carries its weight enumerator. This is the initial step in the suggestion by M. Loebl, to extend the theory of Pfaffian orientations
from graphs to general linear binary codes. Then we carry out also the second step by constructing, for
every triangular configuration $\Delta$, a triangular configuration $\Delta'$ and a bijection between the
cycle space of $\Delta$ and the set of the perfect matchings of $\Delta'$.
\end{abstract}

\end{frontmatter}

\section{Introduction}
A seminal result of Galluccio and Loebl~\cite{loeblpfaffian} asserts that the weight enumerator of the cut space $\mathcal{C}$ of a graph $G$ may be written as a linear combination of $4^{g(G)}$ Pfaffians, where $g(G)$ is the minimal genus of a surface in which $G$ can be embedded.
Recently, a topological interpretation of this result was given by Cimasoni and Reshetikhin~\cite{cimasoni-2007}.

Viewing the cut space $\mathcal{C}$ as a binary linear code, $G$ may be considered as a useful geometric representation of $C$ which provides an important structure for the weight enumerator of $\mathcal{C}$.

This motivated Martin Loebl to ask, about 10 years ago, the following question: Which binary codes are cycle spaces of simplicial complexes?
In general, for the binary codes with a geometric representation, one may hope to obtain a formula analogous to that of Galluccio and Loebl~\cite{loeblpfaffian}.

This question remains open but in this paper we circumvent it by asking for geometric representations which
carry over the weight enumerator only. We note that this weakening is still sufficient for the extension of the theory
of Pfaffian orientations.

We present a construction which shows that a useful geometric representation exists for all binary codes. The main result is as follows:

\begin{thm}
\label{thm:main}
For each binary linear code $\mathcal{C}$ of length $n$, one can construct a triangular configuration $\Delta$ and positive integer $e$ linear in $n$, so that if the weight enumerator of the cycle space of $\Delta$ equals
$$
\sum_{i=0}^{m}a_ix^i
$$
then the weight enumerator of $\mathcal{C}$  satisfies
\begin{equation*}
W_\mathcal{C}(x^2)=\sum_{i=0}^{m}a_ix^{i\bmod e}.
\qed
\end{equation*}
\end{thm}

The second main result of the paper is to construct, for
every triangular configuration $\Delta$, a triangular configuration $\Delta'$ and a bijection between the
cycle space of $\Delta$ and the set of the perfect matchings of $\Delta'$. This carries over the second step in the
Loebl's suggestion to extend the theory of Pfaffian orientations to the general binary linear codes.

%
%

\section{Preliminaries}
\label{sec:prelim}
We will start with definitions of the basic concepts.
Let $n$ be a positive integer. A \dfn{binary linear code $\mathcal{C}$ of length $n$} is a subspace of $GF(2)^n$, and each vector in $\mathcal{C}$ is called a \dfn{codeword}.

The \dfn{weight} of codeword $c$ is the number of nonzero coordinates, denoted by $w(c)$.

A binary linear code $\mathcal{C}$ is \dfn{even} if all codewords have even weight.


We define a partial order on $\mathcal{C}$ as follows:
let $c=(c^1,\dots,c^n),d=(d^1,\dots,d^n)$ be codewords of $\mathcal{C}$. Then $c\preceq d$ if $c^i=1$ implies $d^i=1$ for all $i=1,\dots,n$. Codeword $d$ is \dfn{minimal} if $c\preceq d$ implies $c=d$ for all.

The \dfn{weight enumerator} of the code $\mathcal{C}$ is defined according to the formula
\begin{equation*}
W_\mathcal{C}(x):=\sum_{c\in \mathcal{C}} x^{w(c)}.
\end{equation*}

An \dfn{abstract simplicial complex} on a finite set $V$ is a hereditary family $\Delta$ of subsets of $V$.


Let $X$ be an element of $\Delta$. The \dfn{dimension} of $X$ is $\left\lvert X\right\rvert-1$, denoted by $\dim X$.
The dimension of $\Delta$ is $\max\left\lbrace\dim X \vert X\in\Delta\right\rbrace$, denoted by $\dim \Delta$.


A \dfn{simplex} in $\mathbb{R}^n$ is the convex hull of an affine independent set $V$ in $\mathbb{R}^d$. The \dfn{dimension} of the simplex is $\left\lvert V\right\rvert-1$.
The convex hull of any nonempty subset of $V$ that define a simplex is called a \dfn{face} of the simplex.
A \dfn{simplicial complex} $\Delta$ is a set of simplices fulfilling the following conditions:
\begin{itemize}
 \item Every face of a simplex from $\Delta$ belongs to $\Delta$.
 \item The intersection of every two simplices of $\Delta$ is a face of both.
\end{itemize}

We denote the subset of $d$-dimensional simplices of $\Delta$ by $\Delta^d$.

Every simplicial complex defines an abstract simplicial complex on the set of vertices $V$, namely the family of sets of vertices of simplexes of $\Delta$. We will denote this abstract simplicial complex by $\mathcal{A}(\Delta)$.

The \dfn{geometric realization} of an abstract simplicial complex $\Delta$ is a simplicial complex $\Delta'$ such that $\Delta=\mathcal{A}(\Delta')$.
It is well known that every finite $d$-dimensional abstract simplicial complex can be realized as a simplicial complex in $\mathbb{R}^{2d+1}$. We will denote this geometric realization of an abstract simplicial complex $\Delta$ by $\mathcal{G}(\Delta)$.

This paper studies 2-dimensional simplicial complexes where each maximal simplex is a triangle. We will call them \dfn{triangular configurations}.

The number of triangles in an (abstract) simplicial complex $\Delta$ will be denoted by $\left\lvert\Delta\right\rvert$.

A \dfn{subconfiguration} of a triangular configuration $\Delta$ is a triangular configuration $\Delta'$ such that $\Delta'\subseteq\Delta$.

A \dfn{cycle} of a triangular configuration is a subconfiguration such that every edge is incident with an even number of triangles. A \dfn{circuit} is a minimal non-empty cycle under inclusion.

Let $\Delta_1$, $\Delta_2$ be subconfigurations of a triangular configuration $\Delta$.
The \dfn{difference} of $\Delta_1$ and $\Delta_2$, denoted by $\Delta_1-\Delta_2$, is defined to be the triangular configuration obtained from $\Delta_1^0\cup\Delta_1^1\cup\Delta_1^2\setminus\Delta_2^2$ by removing the edges and vertices that are not contained in any triangle in $\Delta_1^2\setminus\Delta_2^2$.

The \dfn{symmetric difference} of $\Delta_1$ and $\Delta_2$, denoted by $\Delta_1\bigtriangleup\Delta_2$, is defined to be $\Delta_1\bigtriangleup\Delta_2:=\left(\Delta_1\cup\Delta_2\right)-\left(\Delta_1\cap\Delta_2\right)$.

Let $\Delta_1,\Delta_2$ be triangular configurations. The union of $\Delta_1,\Delta_2$ is defined to be $\Delta_1\cup\Delta_2:=\mathcal{G}(\mathcal{A}(\Delta_1)\cup\mathcal{A}(\Delta_1))$.


Let $\Delta$ be a $d$-dimensional simplicial complex. We define the \dfn{incidence matrix} $A=\left(A_{ij}\right)$ as follows: the rows are indexed by $\left(d-1\right)$-dimensional simplices and the columns by $d$-dimensional simplices. We set
\begin{equation*}
a_{ij}:=\begin{cases}
         1& \text{if }(d-1)\text{-simplex }i\text{ belongs to }d\text{-simplex }j,\\
         0& \text{otherwise}.
        \end{cases}
\end{equation*}

The \dfn{cycle space} $\mathcal{C}$ of $\Delta$ is the kernel $\ker\Delta$ of the incidence matrix of $\Delta$ over $GF(2)$, and $\mathcal{C}=\ker\Delta$ is said to be \dfn{represented} by $\Delta$.

For a subconfiguration $C$ of $\Delta$, we let $\chi(C)=(\chi(C)^{t_1},\dots,\chi(C)^{t_{\lvert\Delta\rvert}})\in\left\lbrace 0,1\right\rbrace^{\left\lvert\Delta\right\rvert}$ denote its \dfn{incidence vector}, where $\chi(C)^t:=1$ if $C$ contains triangle $t$, and $\chi(C)^t:=0$ otherwise.


It is well known that the kernel of $\Delta$ is the set of incidence vectors of cycles of $\Delta$.

Let $\mathcal{C}\subseteq\left\lbrace 0,1\right\rbrace^n$ be a binary linear code. Let $S$ be a subset of $\left\lbrace 1,\dots,n\right\rbrace $. \dfn{Puncturing} code $\mathcal{C}$ along $S$ means deleting the entries indexed by the elements of $S$ from each codeword of $\mathcal{C}$. The resulting code is denoted by $\mathcal{C}/S$.


%

\section{Triangular representation of binary codes}
\label{sec:representation}
First, we define three useful triangular configurations.
\subsection{Triangular configuration $B^n$}
\begin{figure}[h]
\begin{center}
 \includegraphics[width=150pt]{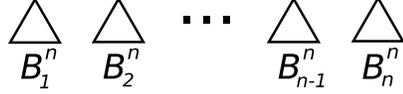}
\end{center}
\caption{Triangular configuration $B^n$.}
\label{fig:trianconf B}
\end{figure}

The triangular configuration $B^n$ consists of $n$ disjoint triangles as is depicted in Figure~\ref{fig:trianconf B}. We denote the triangles of $B^n$ by $B^n_1,\dots,B_n^n$.

\subsection{Triangular sphere $\mathcal{S}^m$}

\begin{figure}[h]
\begin{center}
	\includegraphics[width=50pt]{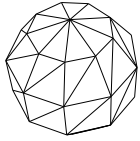}
\end{center}
\caption{Triangular sphere $\mathcal{S}^m$.}
\label{fig:sphere}
\end{figure}

The triangular sphere $\mathcal{S}^m$, depicted in Figure~\ref{fig:sphere}, is a triangulation of a $2$-dimensional sphere by $m$ triangles. This triangulation exists for every even $m\geq 4$.
We denote the triangles of $\mathcal{S}^m$ by $\mathcal{S}^m_1,\dots,\mathcal{S}^m_m$

\subsection{Triangular tunnel $T$}

\begin{figure}[h]
	\begin{center}
	\includegraphics[width=320pt]{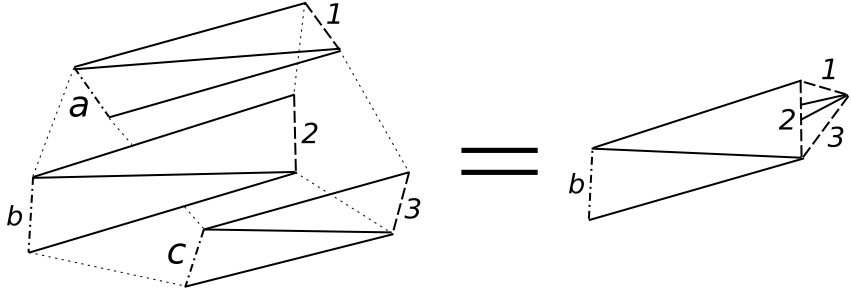}
	\end{center}
\caption{Triangular tunnel $T$.}
\label{fig:tunnel}
\end{figure}

The triangular tunnel $T$ is depicted in Figure~\ref{fig:tunnel}.

In particular, triangles $\left\lbrace 1,2,3\right\rbrace $ and $\left\lbrace a,b,c\right\rbrace $ are not elements of $T$.


\subsection{Joining triangles by tunnels}

Let $\Delta$ be a triangular configuration. Let $t_1$ and $t_2\in\Delta$ be two disjoint triangles of $\Delta$. The join of $t_1$ and $t_2$ in $\Delta$ is the triangular configuration $\Delta'$ defined as follows.
Let $T$ be a triangular tunnel as in Figure~\ref{fig:tunnel}. Let $t_1^1,t_1^2,t_1^3$ and $t_2^1,t_2^2,t_2^3$ be edges of $t_1$ and $t_2$, respectively.
We relabel edges of $T$ such that $\left\{a,b,c\right\}=\left\{t_1^1,t_1^2,t_1^3\right\}$ and $\left\{1,2,3\right\}=\left\{t_2^1,t_2^2,t_2^3\right\}$. Then $\Delta'$ is defined to be $\Delta\cup T$.

\subsection{Construction}

Let $\mathcal{C}$ be a binary code of length $n$ and dimension $d$.
Let $B=\left\{b_1,\dots,b_d\right\}$ be a basis of $\mathcal{C}$.
We construct its triangular representation $\Delta^C_B$ as follows.

\begin{figure}[hp]
\begin{center}
 \includegraphics[width=150pt]{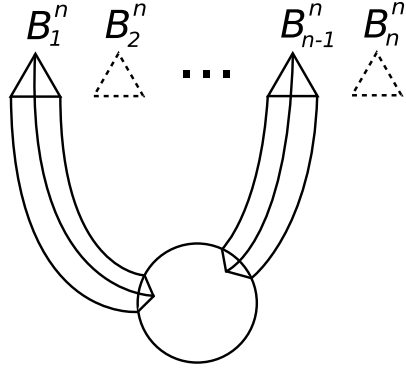}
\end{center}
\caption{$\Delta^\mathcal{C}_{b_i}$ represents a basis vector $\left(1,0,\dots,1,0\right)$ of $\mathcal{C}$.}
\label{fig:trianrep cycle}
\end{figure}

For every basis vector $b_i$ we construct a triangular configuration $\Delta^C_{b_i}$.

\begin{figure}[hp]
\begin{center}
 \includegraphics[width=150pt]{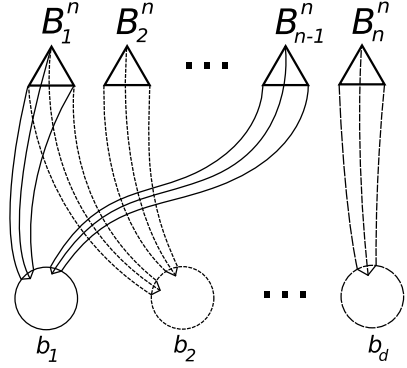}
\end{center}
\caption{An example of triangular representation $\Delta^\mathcal{C}_B$ of $\mathcal{C}$.}
\label{fig:trianrep}
\end{figure}

Triangular configuration $\Delta^\mathcal{C}_{b_i}$ is obtained from $B^n\cup \mathcal{S}^m$, where $m$ is even and $m\geq n$, $m\geq 4$. 
Let $J^i$ be the set of indices of nonzero entries of $b_i$. For each $j\in J^i$ we join triangle $\mathcal{S}^m_j$ of $\mathcal{S}^m$ with triangle $B^n_j$. Then, we remove triangle $\mathcal{S}^m_j$ from $\mathcal{S}^m$.
Finally, we remove triangles of $B^n$ that are not joined with the sphere.
The example of $\Delta^\mathcal{C}_{b_i}$ for $b_i=\left(1,0,\dots,1,0\right)$ is depicted in Figure~\ref{fig:trianrep cycle}.

Thus triangular configuration $\Delta^\mathcal{C}_{b_i}$ contains $B^n_j$ if and only if $j\in J^i$. We note that

\begin{prop}
\label{o:eve}
$|\Delta^\mathcal{C}_{b_i}|-w(b_i)$ is always even.
\end{prop}

Triangular configurations $\Delta^\mathcal{C}_{b_i}$, $i=1,\dots,d$, share the triangles of $B^n$ and do not share spheres $\mathcal{S}^m$. Hence, $\mathcal{A}(\Delta^\mathcal{C}_{b_i})\cap\mathcal{A}(\Delta^\mathcal{C}_{b_j})\subseteq\mathcal{A}(B_n)$ holds for $i<j$, $i,j\in\left\{1,\dots,d\right\}$.

Finally, triangular representation $\Delta^\mathcal{C}_B$ of $\mathcal{C}$ is the union of $\Delta^\mathcal{C}_{b_i}$, $i=1,\dots,d$.

An example of triangular representation $\Delta^\mathcal{C}_B$ of $\mathcal{C}$ is depicted in Figure~\ref{fig:trianrep}.

Triangular representation $\Delta^\mathcal{C}_B$ of $\mathcal{C}$ is \dfn{balanced} if there is an integer $e$ such that $\left\lvert\Delta^\mathcal{C}_{b_i}\right\rvert-w(b_i)=e$ for all $i=1,\dots,d$.
This $e$ is denoted by $\ex{\Delta^\mathcal{C}_B}$.


Let $c$ be a codeword of $\mathcal{C}$. Let $c=\sideset{}{^2}\sum_{i\in I}b_i$ be the unique expression of $c$, where $b_i\in B$. The \dfn{degree} of $c$ with respect to basis $B$ is defined to be the cardinality $\left\rvert I\right\lvert$ of the index set.
The degree is denoted by $\vd{c}$.

We denote by $\ker\Delta^\mathcal{C}_B$ the cycle space of the triangular configuration
$\Delta^\mathcal{C}_B$
We define linear mapping $f\colon\mathcal{C}\mapsto \ker\Delta^\mathcal{C}_B$ in the following way:
let $c$ be a codeword of $\mathcal{C}$; let $c=\sideset{}{^2}\sum_{i\in I}b_i$ be the unique expression of $c$, where $b_i\in B$.
We define $f(c):=\chi(\bigtriangleup_{i\in I} \Delta^\mathcal{C}_{b_i})$.
Hence, the entries of $f(c)$ are indexed by the triangles of $\Delta^\mathcal{C}_B$.
We have $f(c)^{B^n_j}=1$ if and only if $\bigtriangleup_{i\in I} \Delta^\mathcal{C}_{b_i}$ contains triangle $B^n_j$.

\begin{prop}
\label{prop:mappingf}
Denote $\left\lvert \bigtriangleup_{i\in I} \Delta^\mathcal{C}_{b_i} \right\rvert$ by $m$.
Let $c=(c^1,\dots,c^n)$ and
\begin{equation*}
f(c)=\left(f(c)^{B^n_1},\dots,f(c)^{B^n_n},f(c)^{n+1},\dots,f(c)^{m}\right). 
\end{equation*}
Then, $f(c)^{B^n_j}=c^j$ for all $j=1,\dots,n$ and all $c\in\mathcal{C}$.
\end{prop}
\begin{pf}
We show the proposition by induction on the degree $\vd{c}$ of $c$.
Codeword $c$ is equal to $\sideset{}{^2}\sum_{i\in I}b_i$.

If $\vd{c}=0$, then $c=0$ and $f(c)=0$ is the incidence vector of the empty triangular configuration. So, the proposition holds for vectors of degree $0$.

If $\vd{c}$ is greater than $0$, then $\lvert I\rvert\geq 1$. We choose some $k$ from $I$. Codeword $c+^2b_k$ has degree less than $c$. By induction assumption, the proposition holds for $c+^2b_k$.
Let $b_k=(b_k^1,\dots,b_k^n)$.

From the definition of $\Delta^\mathcal{C}_{b_k}$, $b_k^j=\chi(\Delta^\mathcal{C}_{b_k})^{B^n_j}$ for all $j=1,\dots,n$.

Therefore,
\begin{equation*}
c^j=(c^j+^2b_k^j)+^2b_k^j=\chi(\bigtriangleup_{i\in I\setminus\{k\}}\Delta^\mathcal{C}_{b_i})^{B^n_j}+^2\chi(\Delta^\mathcal{C}_{b_k})^{B^n_j}=f(c)^{B^n_j}
\end{equation*}
for all $j=1,\dots,n$. \qed
\end{pf}

\begin{cor}
\label{cor:inj}
Mapping $f$ is injective.
\end{cor}

\begin{lem}
\label{lem:cycletriangles}
Every non-empty cycle of $\Delta^\mathcal{C}_B$ contains $\Delta^\mathcal{C}_{b_i}-B^n$ as a subconfiguration for some $i\in \left\{1,\dots,d\right\}$.
\end{lem}
\begin{pf}
Every cycle of $\Delta^\mathcal{C}_B$ contains either all triangles or no triangle of $\Delta^\mathcal{C}_{b_i}-B^n$, since $\Delta^\mathcal{C}_{b_i}\cap \Delta^\mathcal{C}_{b_j}\subseteq B^n$ for all distinct $i,j\in \left\{1,\dots,d\right\}$.

$B^n$ does not contain non-empty cycles, since the triangles of $B^n$ are disjoint.
Therefore, every non-empty cycle contains a triangle of $\Delta^\mathcal{C}_{b_i}- B^n$ for some $i\in \left\{1,\dots,d\right\}$. Hence, every non-empty cycle contains $\Delta^\mathcal{C}_{b_i}- B^n$. \qed
\end{pf}


\begin{thm}
\label{thm:trianrep}
Let $\mathcal{C}$ be a binary code. Let $\Delta^\mathcal{C}_B$ be its triangular representation with respect to a basis $B$.
Mapping $f$ defined above is a bijection of the binary linear codes 
$\mathcal{C}$ and $\ker\Delta^\mathcal{C}_B$ which maps minimal codewords to
minimal codewords.
\end{thm}
\begin{pf}
By Corollary~\ref{cor:inj}, $f$ is injective.


It remains to be proven that $\dim\mathcal{C}=\dim\ker\Delta^\mathcal{C}_B$.

Suppose on the contrary that some codeword of $\ker\Delta^\mathcal{C}_B$ is not in the span of $\left\{f(b_1),\dots,f(b_d)\right\}$.
Let $c$ be such a codeword with the minimal possible weight $w(c)$.
Let $K$ be a cycle of $\Delta^\mathcal{C}_B$ such that $\chi(K)=c$.

By Lemma~\ref{lem:cycletriangles}, cycle $K$ contains $\Delta^\mathcal{C}_{b_i}- B^n$ for some $i\in \left\{1,\dots,d\right\}$.
Since $\left\lvert\Delta^\mathcal{C}_{b_i}- B^n\right\rvert>\left\lvert B^n\right\rvert$, $\left\lvert K\bigtriangleup \Delta^\mathcal{C}_{b_i}\right\rvert < \left\lvert K\right\rvert$.
Therefore, $w(c)>w(\chi(K\bigtriangleup \Delta^\mathcal{C}_{b_i}))$. This is a contradiction.


Finally we show that $f$ maps minimal codewords to minimal codewords.
Let $d$ be a minimal codeword.
Suppose on the contrary that $f(d)$ is not a minimal codeword of $\ker\Delta^\mathcal{C}_B$.
Then for some codeword $c$, $f(c)\prec f(d)$. However, $c^i=f(c)^i=1$ implies that $d^i=f(d)^i=1$, so $c\prec d$. This contradicts the minimality of $d$. \qed

\end{pf}

\begin{figure}[hp]
 \centering
 \includegraphics[width=150pt]{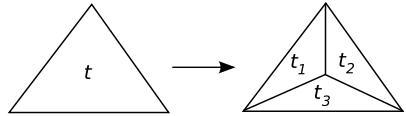}
 \caption{Triangle subdivision}
 \label{fig:trian div}
\end{figure}

Let $t$ be a triangle of a triangular configuration $\Delta$. A \dfn{subdivision} of triangle $t$ is a triangular configuration obtained from $\Delta$ by exchanging triangle $t$ is replaced by triangles $t_1,t_2,t_3$ in the way depicted in Figure~\ref{fig:trian div}.

\begin{prop}
\label{prop:balanced}
Every binary code $\mathcal{C}$ of length $n$ and dimension $d$ has a balanced triangular representation $\Delta^\mathcal{C}_B$ such that $\ex{\Delta^\mathcal{C}_B}>n$, where $B$ is an arbitrary basis for $\mathcal{C}$.
\end{prop}
\begin{pf}
Let $\Delta^\mathcal{C}_B$ be an arbitrary triangular representation of $\mathcal{C}$ with respect to a basis $B=\left\{b_1,\dots,b_d\right\}$.
We denote by $k_i$ the number $\left\lvert\Delta^\mathcal{C}_{b_i}\right\rvert-w(b_i)$.
Every $k_i$ is even by Proposition~\ref{o:eve}.
Let $n'$ be the smallest even number greater than $n$.
Let $k:=\max\left\{n',k_i \vert i=1,\dots,d\right\}$.

For each $i\in\left\{1,\dots,d\right\}$ such that $k_i\neq k$, the following step is applied.
We choose a triangle $t$ from $\Delta^\mathcal{C}_{b_i}- B^n$ and subdivide $t$.
The number $k_i$ is increased by $2$.
If $k_i$ still does not equal to $k$, then we repeat this step.

After this procedure, $\Delta^\mathcal{C}_B$ is balanced and $\ex{\Delta^\mathcal{C}_B}>n$. \qed
\end{pf}


\begin{prop}
\label{prop:weightf}
Let $\mathcal{C}$ be an even binary linear code. Let $\Delta^\mathcal{C}_B$ be its balanced triangular representation with respect to a basis $B$. Then $w(f(c))=w(c)+\vd{c}\ex{\Delta^\mathcal{C}_B}$ for every codeword $c\in \mathcal{C}$.
\end{prop}
\begin{pf}
Write $c$ as $\sideset{}{^2}\sum_{i\in I}b_i$, where $b_i\in B$.
Then $f(c)=\chi(\bigtriangleup_{i\in I} \Delta^\mathcal{C}_{b_i})$.
Now, $\bigtriangleup_{i\in I} \Delta^\mathcal{C}_{b_i}$ contains all triangles of $\Delta^\mathcal{C}_{b_i}- B^n$ for all $i\in I$.
The number of these triangles is $\vd{c}\ex{\Delta^\mathcal{C}_B}$, since $\left\lvert\Delta^\mathcal{C}_{b_i}- B^n\right\rvert=\ex{\Delta^\mathcal{C}_B}$ and $\left\lvert I\right\rvert=\vd{c}$.

By Proposition~\ref{prop:mappingf}, $\bigtriangleup_{i\in I} \Delta^\mathcal{C}_{b_i}$ contains triangle $B^n_k$ if and only if $c_k=1$. The number of these triangles is $w(c)$.

Therefore, $w(f(c))=w(c)+\vd{c}\ex{\Delta^\mathcal{C}_B}$. \qed
\end{pf}

\section{Weight enumerator}
In this section, we state the connection between the weight enumerator of a code and the weight enumerator of its triangular representation. This will provide a proof of Theorem~\ref{thm:main}.


We define the \dfn{extended weight enumerator} (with respect to a fixed basis) by
\begin{equation*}
W^k_\mathcal{C}(x):=\sum_{\substack{c\in \mathcal{C}\\\vd{c}=k}} x^{w(c)}.
\end{equation*}

If code $\mathcal{C}$ has dimension $d$, then
\begin{equation*}
W_\mathcal{C}(x)=\sum_{k=0}^d W^k_\mathcal{C}(x).
\end{equation*}

\begin{prop}
\label{prop:extended polynomial}
Let $\mathcal{C}$ be a binary code and let $\Delta^\mathcal{C}_B$ be its balanced triangular representation $\Delta^\mathcal{C}_B$ with respect to the fixed basis $B$. Then
\begin{equation*}
W^k_{\ker\Delta^\mathcal{C}_B}(x)=W^k_\mathcal{C}(x)x^{k\ex{\Delta^\mathcal{C}_B}}.
\end{equation*}
\end{prop}
\begin{pf}
Let $f$ be the mapping defined in Section~\ref{sec:representation}.
For every codeword $c$ of degree $k$ of $\mathcal{C}$ there is codeword $f(c)$ of degree $k$ of $\ker\Delta^\mathcal{C}_B$.
By Proposition~\ref{prop:weightf}, $w(f(c))=w(c)+k\ex{\Delta^\mathcal{C}_B}$.

Therefore,
\begin{equation*}
\begin{split}
W^k_{\ker\Delta^\mathcal{C}_B}(x)& =\sum_{\substack{f(c)\in \ker\Delta^\mathcal{C}_B\\\vd{f(c)}=k}} x^{w(f(c))} =\sum_{\substack{c\in \mathcal{C}\\\vd{c}=k}} x^{w(c)+k\ex{\Delta^\mathcal{C}_B}}
 =W^k_{\mathcal{C}}(x)x^{k\ex{\Delta^\mathcal{C}_B}}.
\qed
\end{split}
\end{equation*}
\end{pf}
\begin{prop}
\label{prop:weights}
Let $\mathcal{C}$ be a binary code of length $n$. Let $\Delta^\mathcal{C}_B$ be a balanced triangular representation of $\mathcal{C}$.
The inequality $k\ex{\Delta^\mathcal{C}_B}\leq w(c)\leq k\ex{\Delta^\mathcal{C}_B}+n$ holds for every codeword $c$ of degree $k$ of $\ker\Delta^\mathcal{C}_B$.
\end{prop}
\begin{pf}
By Proposition~\ref{prop:weightf}, $w(c)=w(f^{-1}(c))+k\ex{\Delta^\mathcal{C}_B}$.
Since $0\leq w(f^{-1}(c))\leq n$ holds for every $c\in \ker\Delta^\mathcal{C}_B$, $k\ex{\Delta^\mathcal{C}_B}\leq w(c)\leq k\ex{\Delta^\mathcal{C}_B}+n$. \qed
\end{pf}

\begin{cor}
Let $\mathcal{C}$ be a binary code of dimension $d$ and length $n$. Let $\Delta^\mathcal{C}_B$ be a balanced triangular representation of $\mathcal{C}$ such that $n<\ex{\Delta^\mathcal{C}_B}$. Denote $\ex{\Delta^\mathcal{C}_B}$ by $e$. If $\sum_{i=0}^{de+n}a_ix^i$ is the weight enumerator of $\ker\Delta^\mathcal{C}_B$, then
\begin{equation*}
W^k_{\ker\Delta^\mathcal{C}_B}(x)=\sum_{i=ke}^{ke+n}a_ix^i.
\end{equation*}
\end{cor}
\begin{pf}
By Proposition~\ref{prop:weights}, $w(c)\leq(k-1)e+n$ holds for all codewords $c\in \ker\Delta^\mathcal{C}_B$ of degree less than $k$. Since $n<e$, $w(c)\leq ke-e+n<ke$.

By Proposition~\ref{prop:weights}, $(j+1)e\leq w(c)$ holds for all codewords $c\in \ker\Delta^\mathcal{C}_B$ of degree greater than $k$. Since $n<e$, $ke+e<ke+n\leq w(c)$.

Hence, enumerator $W^k_{\ker\Delta^\mathcal{C}_B}(x)$ is the sum over all codewords of weight between $ke$ and $ke+n$. \qed
\end{pf}

\begin{thm}
\label{thm:weight polynomial}
Let $\mathcal{C}$ be a binary code of dimension $d$ and length $n$. Let $\Delta^\mathcal{C}_B$ be a balanced triangular representation of $\mathcal{C}$ such that $n<\ex{\Delta^\mathcal{C}_B}$. Denote $\ex{\Delta^\mathcal{C}_B}$ by $e$. If $\sum_{i=0}^{de+n}a_ix^i$ is the weight polynomial of $\ker\Delta^\mathcal{C}_B$, then
\begin{equation*}
W_\mathcal{C}(x)=\sum_{i=0}^{de+n}a_ix^{i\bmod e}.
\end{equation*}
\end{thm}
\begin{pf}
$w(c)\leq n$ for every codeword $c\in \mathcal{C}$.

Let $f$ be the mapping defined in Section~\ref{sec:representation}.
By Proposition~\ref{prop:weightf}, $w(f(c))=w(c)+\vd{c}e$ for every codeword $c$ of $\mathcal{C}$.
Since $n<e$,
\begin{equation*}
w(f(c))\bmod e = (w(c)+\vd{c}e)\bmod e = w(c).
\end{equation*}

Hence,
\begin{equation*}
W_\mathcal{C}(x)=\sum_{i=0}^{de+n}a_ix^{i\bmod e}.
\qed
\end{equation*}
\end{pf}

Finally, we prove Theorem~\ref{thm:main}.
\begin{pf*}{PROOF of Theorem~\ref{thm:main}.}
Let $\mathcal{C}$ be a linear binary code of length $n$. 
 By Proposition~\ref{prop:balanced}, we can construct a balanced triangular representation $\Delta$ of $\mathcal{C}$ such that $\ex{\Delta}>n$.

Denote $\ex{\Delta}$ by $e$. Let $W_{\Delta}(x)=\sum_{i=0}^{de+n}a_ix^i$ be the weight enumerator of $\Delta$.


Hence,
\begin{equation*}
W_\mathcal{C}(x)=\sum_{i=0}^{de+n}a_ix^{i\bmod e}.
\qed
\end{equation*}
\end{pf*}

\section{Matching}

In this section we reduce computation of the weight enumerator of the even subconfigurations to computation of the weight enumerator of the perfect matchings.

Let $\Delta$ be a triangular configuration. A \dfn{matching} of $\Delta$ is a subconfiguration $M$ of $\Delta$ such that $t_1\cap t_2$ does not contain an edge for every distinct $t_1,t_2\in T(M)$.

Let $\Delta$ be a triangular configuration. Let $M$ be a matching of $\Delta$. Then the \dfn{defect} of $M$ is the set $E(T)\setminus E(M)$. We denote matching with this defect by $M_{E(T)\setminus E(M)}$.

The \dfn{perfect matching} of $\Delta$ is a matching with empty defect. We denote the set of all perfect matchings of $\Delta$ by $\mathcal{P}(\Delta)$.


The weight enumerator of perfect matchings in $\Delta$ is defined to be $P_\Delta(x)=\sum_{P\in\mathcal{P}(\Delta)} x^{w(P)}$, where $w(P):=\sum_{t\in P} w_t$.

%
%
%
%

\subsection{Triangular configuration $P$}
\begin{figure}[ht]
 \centering
 \includegraphics{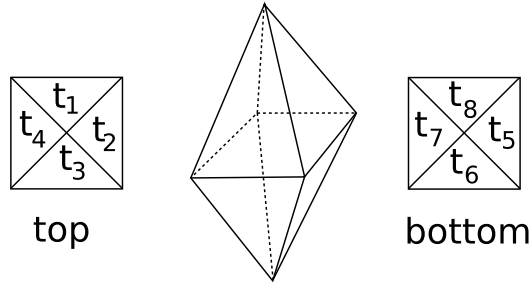}
 \caption{Triangular configuration $P$}
 \label{fig:piramid}
\end{figure}
The triangular configuration $P$ is depicted in Figure~\ref{fig:piramid}.
\begin{prop}
Triangular configuration $P$ has exactly two perfect matchings $\{t_1,t_3,t_5,t_7\}$, $\{t_2,t_4,t_6,t_8\}$.
\end{prop}
\subsection{Closed triangular tunnel $T$}

\begin{figure}[h]
	\begin{center}
	\includegraphics[width=320pt]{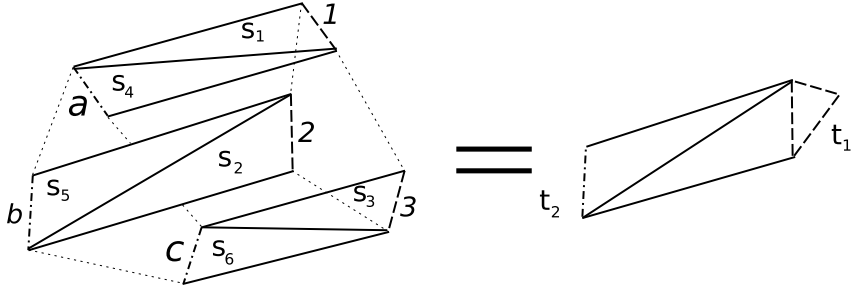}
	\end{center}
\caption{Closed triangular tunnel $T$.}
\label{fig:tunnel2}
\end{figure}

Closed triangular tunnel $T$ is depicted in Figure~\ref{fig:tunnel2}. We call triangles $\{a,b,c\}=t_2$ and $\{1,2,3\}=t_1$ ending triangles.

\begin{prop}
Closed triangular tunnel $T$ has two perfect matchings $M^T_{t_1}=\{t_1,s_4,s_5,s_6\}$, $M^T_{t_2}=\{t_2,s_1,s_2,s_3\}$.
\end{prop}


\subsection{Triangular configuration $E_{pq}$}
\begin{figure}[hp]
 \centering
 \includegraphics{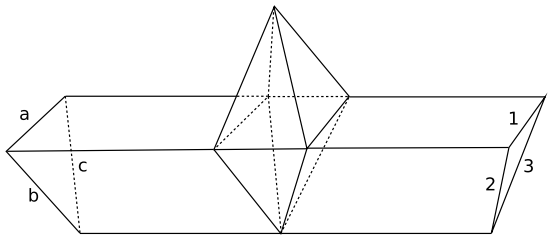}
 \caption{Matching triangular edge}
 \label{fig:matchedge}
\end{figure}

The \dfn{matching triangular edge} is triangular configuration which is obtained from triangular configuration $P$ and two closed triangular tunnels $T$ in the following way:
Let $T_1$ and $T_2$ be closed triangular tunnels and let $t_1^{T_1},p^{T_1}$ and $t_1^{T_2},q^{T_2}$ be ending triangle of $T_1$ and $T_2$, respectively.
We identify $t_1^{T_1}$ with $t_1^P$ and $t_1^{T_2}$ with $t_3^P$. Then $E_{pq}$ is defined to be $T_1\bigtriangleup P\bigtriangleup T_2$.
Triangular configuration $E_{pq}$ is depicted in Figure~\ref{fig:matchedge}.

\begin{prop}
\label{prop:trianedgematch}
Matching triangular edge has two perfect matchings denoted by $N^1_{pq}$ and $N^0_{pq}$.
\end{prop}
\begin{pf}
There are two matchings.
First matching is $N^0_{pq}:=M^{T_1}_{t_1}\cup M^{T_2}_{t_1}\cup \{t_5^P,t_7^P\}$.
Second matching is $N^1_{pq}:=M^{T_1}_{p}\cup M^{T_2}_{q}\cup \{t_2^P,t_4^P,t_6^P,t_8^P\}$.

Any perfect matching of $E_{pq}$ has to contain $\{t_5^P,t_7^P\}$ or $\{t_2^P,t_4^P,t_6^P,t_8^P\}$. This determines remaining triangles in a perfect matching. Hence, there are just two perfect matchings.
\qed
\end{pf}

We denote the matching $N^1_{pq}$ by $M^1_{pq}$ and the matching $N^0_{pq}\setminus{p,q}$ by $M^0_{pq}$.

\subsection{Triangular configuration $T_{pqr}$}
\begin{figure}[hp]
 \centering
 \includegraphics{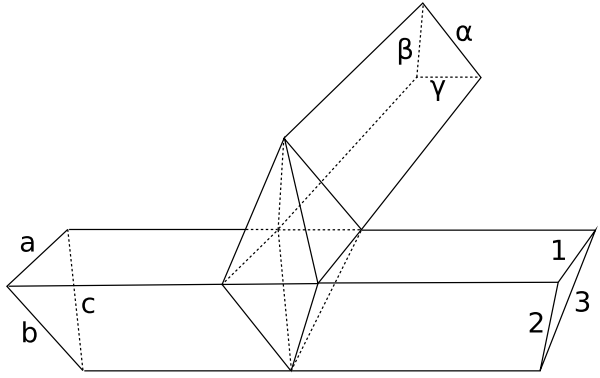}
 \caption{Matching triangle}
 \label{fig:matchtrian}
\end{figure}

The \dfn{matching triangular triangle} is triangular configuration which is obtained from triangular configuration $P$ and three closed triangular tunnels $T$ in the following way:

Let $T_1$, $T_2$ and $T_3$ be closed triangular tunnels and let $t_1^{T_1},p^{T_1}$, $t_1^{T_2},q^{T_2}$ and $t_1^{T_3},r^{T_3}$ be ending triangle of $T_1$, $T_2$ and $T_3$, respectively.
We identify $t_1^{T_1}$ with $t_1^P$, $t_1^{T_2}$ with $t_3^P$ and $t_1^{T_3}$ with $t_5^P$. Then $T_{pqr}$ is defined to be $T_1\bigtriangleup P\bigtriangleup T_2\bigtriangleup T_3$.
Triangular configuration $T_{pqr}$ is depicted in Figure~\ref{fig:matchtrian}.

\begin{prop}
Matching triangular triangle has two perfect matchings denoted by $N^1_{pqr}$ and $N^0_{pqr}$.
\end{prop}
\begin{pf}
There are two matchings.
First matching is $N^0_{pqr}:=M^{T_1}_{t_1}\cup M^{T_2}_{t_1}\cup M^{T_3}_{t_1}\cup \{t_7^P\}$.
Second matching is $N^1_{pqr}:=M^{T_1}_{p}\cup M^{T_2}_{q}\cup M^{T_3}_{r}\cup \{t_2^P,t_4^P,t_6^P,t_8^P\}$.

Any perfect matching of $T_{pqr}$ has to contain $\{t_5^P,t_7^P\}$ or $\{t_2^P,t_4^P,t_6^P,t_8^P\}$. This determines remaining triangles in a perfect matching. Hence, there are just two perfect matchings.
\qed
\end{pf}

We denote the matching $N^1_{pqr}$ by $M^1_{pqr}$ and the matching $N^0_{pqr}\setminus{p,q,r}$ by $M^0_{pqr}$.

\subsection{Triangular configuration $C_{t_1t_2\dots t_n}$}
\begin{figure}[h]
 \centering
 \includegraphics{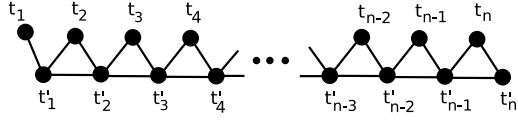}
 \caption{Triangular configuration $C_{t_1t_2\dots t_n}$}
 \label{fig:C}
\end{figure}

This part of reduction is analogous to the reduction for graphs described in Galluccio et al~\cite{galluccio01optimization}.
Let $t^C_1,t'_1$ be empty disjoint triangles.
Let $t_2^C,\dots,t^C_n,t'_2,\dots,t'_n$ be disjoint triangles.
Then $C_{t_1t_2\dots t_n}$ is defined to be $\left( \bigtriangleup_{i=1}^nt_i\right)\bigtriangleup\left(\bigtriangleup_{i=1}^nt'_i\right)\bigtriangleup\left(\bigtriangleup_{i=1}^nE_{t_it'_i}\right)\bigtriangleup\left(\bigtriangleup_{i=2}^nE_{t_it'_{i-1}}\right)\bigtriangleup\left(\bigtriangleup_{i=1}^{n-1}E_{t'_it'_{i+1}}\right)$. The configuration is depicted in Figure~\ref{fig:C}.

\begin{prop}
Let $M^I_C$ denote perfect matching containing triangles $t_i,i\in I$. Then there exists exactly one perfect matching $M^I_C$ of $C_{t_1t_2\dots t_n}$ if and only if $\lvert I\rvert$ is even.
\end{prop}

\begin{pf}
We construct matching by the following algorithm.

In the first step.
If $t_1\in I$ then $M:=M^1_{t_1t'_1}$ else $M:=M^0_{t_1t'_1}$.
We say that first step is even if $t'_1$ is covered otherwise odd. (first step is even if and only if $t_1\notin I$).

In the $i$-th step.
\begin{itemize}
\item If $t'_i$ is covered by $M$ and $t_{i+1}\notin I$ then $M:=M\cup M^0_{t'_it_{i+1}}$.
\item If $t'_i$ is covered by $M$ and $t_{i+1}\in I$ then $M:=M\cup M^0_{t'_it_{i+1}}\cup M^0_{t'_it'_{i+1}}\cup M^1_{t_{i+1}t'_{i+1}}$.
\item  If $t_i$ is not covered by $M$ and $t_{i+1}\notin I$ then $M:=M\cup M^0_{t'_it_{i+1}}\cup M^1_{t'_it'_{i+1}}$.
\item If $t'_i$ is not covered by $M$ and $t_{i+1}\in I$ then $M:=M\cup M^1_{t'_it_{i+1}}$.
\end{itemize}

The $i$-th step is even if $t_i$ is covered. From definition of $i$-th step follows that parity of step changes if $t_{i+1}\in I$.

Hence, the algorithm succed in construction of perfect matching if and only if $\lvert I\rvert$ is even.
\qed
\end{pf}

\subsection{Reduction}

Let $\Delta$ be a triangular configuration.
We construct triangular configuration $\Delta'$ such that every even subconfiguration of $\Delta$ uniqually coresponds to one perfect matching of $\Delta'$ and a natural weight-preserving bijection between the set of the even subconfiguration of $\Delta$ and the set of the perfect matchings of $\Delta'$.

We put into $\Delta'$ empty disjoint triangles $t_e$ for every tuple $(t,e)$ where $e\in E(\Delta)$ and $t\in T(\Delta)$.

We add to $\Delta'$ matching triangle $T_{t_at_bt_c}$ for every triangle $t\in T(\Delta)$, where $a,b,c$ are edges of $t$.
We assign weight $w(t):=1$ to one triangle in matching $M_t^1$ and weight $0$ otherwise.

We add to $\Delta'$ triangular configuration $C_{t^1_e\dots t^n_e}$ for every edge $e\in E(\Delta)$ where $t^1_e,\dots,t^n_e$ are triangles incident with $e$ in $\Delta$.
We assign weight $w_t:=0$ to all triangles of $C_{t^1_e\dots t^n_e}$.

\begin{thm}
Let $\Delta$ be a triangular configuration and $\Delta'$ be matching reduction of $\Delta$.
Let $C$ be an even subconfiguration of $\Delta$. Then there exists exactly one perfect matching $M_C$ in $\Delta'$. In $\Delta'$ are not any others perfect matchings.
\end{thm}

\begin{pf}
Let $C$ be an even subconfiguration of $\Delta$. We construct perfect matchings $M_C$ in $\Delta'$.

We denote matchings $M^1_{t_at_bt_c}$ and $M^0_{t_at_bt_c}$ of $T_{t_at_bt_c}$ by $M_t^1$ and $M_t^0$, respectively.

We denote the set $\{i\vert e\in T(t_i),t_i\in C\}$ by $I_e$.

$M_C:=\{M_t^1\vert t\in C\}\cup\{M_t^0\vert t\notin C, t\in T(\Delta)\}\cup\{M_e^{I_e}\vert e\in E(\Delta)\}$.


Matching $M_C$ is perfect.

We show that there is no other matching.
Every matching triangle $T_t$ has to be covered by $M^1_t$ or $M^0_t$. Thus $C_e$ has to be covered by $M^I_e$ for some even $I$. So every perfect matching in $\Delta'$ defines even subset in $\Delta$.
\qed
\end{pf}

\begin{prop}
Let $\Delta$ be a triangular configuration. Let $\Delta'$ be its matching representation. Let $C$ be an even subconfiguration and $M_C$ be coresponding perfect matching then $\lvert C\rvert=w(M_C)$.
\end{prop}
\begin{pf}
\begin{align*}
 w(M_C)&=\sum_{t\in C}w(M^1_t) + \sum_{t\notin C, t\in T(\Delta)}w(M^0_t) + \sum_{e\in E(\Delta)}w(M_e^{\{i\vert e\in T(t_i),t_i\in C\}})\\
&=\sum_{t\in C}1 + \sum_{t\notin C, t\in T(\Delta)}0 + \sum_{e\in E(\Delta)}0\\
&=\lvert C\rvert
\end{align*}
 \qed
\end{pf}

The following theorem is a consequence.

\begin{thm}
Let $\Delta$ be a triangular configuration. Let $\Delta'$ be its matching representation.
Then $W_\Delta(x)=P_{\Delta'}(x)$.
\end{thm}

\begin{ack}
This article extends a result of my master thesis, written under the direction of Martin Loebl.
I would like to thank Martin Loebl for helpful discussions and continuous support.
\end{ack}

\bibliography{references}{}

\begin{thebibliography}{1}
\expandafter\ifx\csname url\endcsname\relax
  \def\url#1{\texttt{#1}}\fi
\expandafter\ifx\csname urlprefix\endcsname\relax\def\urlprefix{URL }\fi

\bibitem{loeblpfaffian}
A.~Galluccio, M.~Loebl, On the theory of {P}faffian orientations. {I}.
  {P}erfect matchings and permanents., Electronic Journal of Combinatorics 6
  (1999) .

\bibitem{cimasoni-2007}
D.~Cimasoni, N.~Reshetikhin, Dimers on surface graphs and spin structures. i,
  Communications in Mathematical Physics 275~(1) (2007) 187--208.

\bibitem{galluccio01optimization}
A.~Galluccio, M.~Loebl, J.~Vondrák, Optimization via enumeration: a new
  algorithm for the max cut problem, Mathematical Programming 90~(2) (2001)
  273--290.
\newline\urlprefix\url{citeseer.ist.psu.edu/galluccio01optimization.html}

\end{thebibliography}
\bibliographystyle{elsart-num}
\end{document}